\documentclass{amsart}
\usepackage[all]{xy}
\newtheorem{innercustomthm}{Theorem}
\newenvironment{customthm}[1]
  {\renewcommand\theinnercustomthm{#1}\innercustomthm}
  {\endinnercustomthm}

\theoremstyle{plain}
\newtheorem{thm}{Theorem}[section]
\newtheorem{cor}[thm]{Corollary}
\newtheorem{prop}[thm]{Proposition}
\newtheorem{lem}[thm]{Lemma}

\theoremstyle{definition}
\newtheorem{defn}[thm]{Definition}

\newtheorem{rmk}[thm]{Remark}

\newenvironment{customdef}[2][Definition]{\begin{trivlist}
\item[\hskip \labelsep {\bfseries #1}\hskip \labelsep {\bfseries #2}]}{\end{trivlist}}

\def\P{\mathbb{P}}
\def\O{\mathcal{O}}

\begin{document}

\title{Secant spaces and syzygies of special line bundles on curves}
\author{Marian Aprodu}
\email{marian.aprodu@imar.ro}
\address{"Simion Stoilow" Institute of Mathematics of the Romanian Academy, Bucharest}
\author{Edoardo Sernesi}
\email{sernesi@mat.uniroma3.it}
\address{Dipartimento di Matematica e Fisica, Universit\`a Roma Tre, Roma.}
\thanks{The first named author thanks the University of Trento and Dipartimento di Matematica e Fisica, Universit\`a Roma Tre for hospitality during the preparation of this work. The second named author is grateful to the IMAR and to the University of Trento for hospitality during the preparation of this work. The authors have been partly supported by a RIP of CIRM-Trento. The first named author was partly supported by the grant PN-II-ID-PCE-2012-4-015. The second named author was partly supported by the project MIUR-PRIN 2010/11 \emph{Geometria delle variet\`a algebriche}}
\date{\today}

\begin{abstract} On a special line bundle $L$ on a projective curve $C$ we introduce a geometric condition called $(\Delta_q)$. When  $L=K_C$  this condition implies gon$(C) \ge q+2$.  For an arbitrary special $L$  we show that $(\Delta_3)$ implies that $L$ has the well-known property  $(M_3)$, generalizing a  similar result  proved by Voisin in the case $L=K_C$. \end{abstract}

\maketitle

\section{Introduction}

In this paper we introduce some new geometric methods in the study of the Koszul cohomology groups of a projective curve with coefficients in an invertible sheaf.
The basic set-up is as follows.

Let $C$ be a smooth complex projective curve of genus $g$ and $L$  a very ample line bundle of degree $d$ on $C$ with $h^0(C,L)=r+1$. Consider a coherent sheaf $\mathcal{F}$ on $C$ and let $V=H^0(C,L)$;
one has natural complexes of vector spaces:
\[
\xymatrix{
\bigwedge^{p+1}V\otimes H^0(\mathcal{F}\otimes L^{q-1}) \ar[r] &\bigwedge^pV\otimes H^0(\mathcal{F}\otimes L^q) \ar[r] & \bigwedge^{p-1}V\otimes H^0(\mathcal{F}\otimes L^{q+1})}
\]
whose cohomology $K_{p,q}(C,\mathcal{F};L)$ is called   \emph{$(p,q)$ (mixed) Koszul cohomology group} of $C$ with respect to $\mathcal{F}$ and $L$. These vector spaces give information about the minimal resolution of the graded module 
\[
\gamma(C,\mathcal{F};L) = \bigoplus_k H^0(\mathcal{F}\otimes L^k)
\]
over 
  the symmetric (polynomial) algebra $R=S^*V$ in a well known way (see \cite{AN10}). The most important cases are obtained for $\mathcal{F}=\mathcal{O}_C$;   the corresponding   graded
$R$-module $\bigoplus_k H^0(L^k)$ is denoted by $\gamma(C;L)$ and its Koszul cohomology groups by $K_{p,q}(C;L)$.  The choice $L=K_C$ is  of central importance, and its study is at the origin of several results and conjectures on this subject. 
The guiding notions are the so-called properties $(N_p)$.

\begin{defn}
The line bundle $L$ \emph{has property} $(N_0)$ if and only if the natural restriction map $\rho: R\longrightarrow \gamma(C;L)$ is surjective i.e. $L$ is normally generated. 
For $p\ge 1$, we say that the bundle $L$ \emph{satisfies the property} $(N_p)$ if and only if it is normally generated and $K_{i,j}(C;L)=0$ for all $j\ne 1$ and all $1\le i\le p$.
\end{defn}

Roughly speaking, $(N_p)$ holds if and only if the minimal resolution of $\gamma(C;L)$  behaves nicely up to the $p$-th step. These notions have provided an excellent motivation on these problems in two important cases, namely in the case $L=K_C$ and in the case deg$(L) \gg 0$.  As an example we recall the following:

\begin{thm}[Green-Lazarsfeld \cite{GL85}]\label{T:GL1} 
If $\mathrm{deg}(L)\ge 2g+1+p$ then $L$ has property $(N_p)$. If $\mathrm{deg}(L)\ge 2g+p$ then $L$ has property $(N_p)$ unless $C$ is hyperelliptic or $L$ embeds $C$ in $\mathbb P^{g+p}$ with a $(p+2)$--secant $p$--plane.
\end{thm}

On the other hand there is virtually no progress in the study of the properties $(N_p)$ in the case when $L$ is a \emph{special} line bundle different from $K_C$. Already $(N_0)$ and $(N_1)$ have escaped any systematic classification for obvious reasons: normal generation and ideal generation of special projective curves behave essentially wildly and it is therefore very difficult to get even a conjectural picture of how the resolution of $\gamma(C;L)$ might look like (see \cite{AN10} section 4.4 for a short discussion).

However a possible solution comes from the study of other properties of $\gamma(C;L)$,  called $(M_q)$, that have been introduced in \cite{GL86} for $q \ge 1$. We shall work with a slightly weaker condition than in loc.cit., in the spirit of \cite{sE94}.

\begin{customdef}{\ref{def:Mq}}
The line bundle $L$ \emph{has property} $(M_q)$ if $K_{n,1}(C;L)=0$ for all $n\ge r-q$.
\end{customdef}

These  are properties enjoyed by the \emph{tail} of the resolution of $\gamma(C;L)$; i.e. property $(M_q)$ holds for $L$ if the resolution of $\gamma(C;L)$ has a nice behaviour at the last $q$ steps. Another, perhaps more suggestive, point of view consists in considering the resolution of the module $\gamma(C,K_C;L)$. Since it is dual to $\gamma(C;L)$ properties $(M_q)$ for $L$ have to do with nice behaviour of the \emph{head} of the resolution of $\gamma(C,K_C;L)$. In his landmark paper \cite{kP25} Petri had already pointed his attention on the module $\gamma(C,K_C;L)$ \emph{when $L$ is special}. In \cite{AS78} the authors showed that Petri's analysis contains a proof of $(M_1)$ for all $L$ on a non--rational curve $C$ and a characterisation of the validity of $(M_2)$ when $L$ is special. Note that when $L=K_C$,  the self-duality of the resolution of $\gamma(C;K_C)$ implies that property $(M_q)$ is equivalent to property $(N_{q-1})$, so that the result discussed in \cite{AS78} generalises Petri's celebrated analysis of the ideal of the canonical model of a non-hyperelliptic curve (see \cite{SD73}). 

The present paper is devoted to the study of $(M_3)$ for a special $L$. This property has been already studied and characterised for $L=K_C$ by Schreyer in \cite{Sch91}, by Voisin in \cite{cV88} and when $\deg(L) \gg 0$ by Ehbauer in \cite{sE94}.  The main issue in considering the  case of any special line bundle, not considered by them,  is to find natural geometric conditions on $C$ and $L$. We introduce the following definition.

\begin{customdef}{\ref{D:deltaq}}
Assume that $ r\ge 4$ and let $1\le q \le 1+r/2$. We say that a reduced effective divisor 
$D = x_1+\cdots + x_{r-q+2}$ on $C$ satisfies condition $(\Delta_q)$ with respect to $L$ if the following conditions are satisfied:

\begin{itemize}
	\item[(a)] $h^0(L(-D)) = q$, 
	\item[(b)] $L(-D)$ is base point free,
	\item[(c)] $h^0(L(-D+x_i))=h^0(L(-D))$ for all $i=1,\dots, r-q+2$.
\end{itemize}
\end{customdef}

In the case $L=K_C$ a divisor $D$ satisfies condition $(\Delta_q)$ if it defines a primitive $g^1_{g-q+1}$. In general $D$ defines an  $(r-q)$-plane in $\P^r$ which is precisely $(r-q+2)$-secant to $\varphi_L(C) \subset \P^r$.
This condition has appeared   in \cite{GL85} in the case $q=2$ and in \cite{cV88}, where it is called $(H_1)$, in the case $q=3$. In both cases they have proved to be the key for $(M_2)$ to hold for $K_C$ (equivalent to Petri's theorem) and $(M_3)$ to hold for $K_C$ respectively. More precisely, a divisor $D=x_1+ \cdots + x_{g-1}$ satisfying condition $(\Delta_2)$ for $K_C$ defines a primitive $g^1_{g-1}$ and the existence of such a $D$ can be seen to be   equivalent to $C$ being not exceptional, i.e. to Cliff$(C)\ge 2$: this is how Green and Lazarsfeld arrive to Petri's Theorem involving $(\Delta_2)$ and using Mumford-Martens' Theorem. On the other hand Voisin shows, for $g\ge 11$ via an elaborate analysis, that $(\Delta_3)$ plus Cliff$(C) \ge 3$ imply that a general projection in $\P^5$ of the canonical model of $C$ satisfies $(M_3)$. It is interesting to note that this is achieved  by excluding in particular that the projected curve lies in certain surfaces that are intersection of quadrics in $\P^5$. Here one cannot but observe the analogy with the way Ehbauer proved $(M_3)$ for $L$ such deg$(L) \gg 0$ in \cite{sE94}: while his method is  different from Voisin's, he is led to consider the same list of surfaces.

Our main result involves condition $(\Delta_3)$ plus a transversality condition as the key hypothesis. Specifically we prove the following:

\begin{customthm}{\ref{T:m3}}
Assume $g\ge 14$, $r \ge 5$, that $L$ is very ample and special of degree $\ge r+13$,    that    each component of the locus of $(r-1)$-secant $(r-3)$-planes has the expected dimension $r-4$, and that the general such $(r-3)$-plane  in each component satisfies $(\Delta_3)$
with respect to $L$. Then $L$ satisfies $(M_3)$ unless $\mathrm{Cliff}(C) \le 2$.
\end{customthm}

 The relation between Condition $(\Delta_3)$ and the vanishing of the $K_{i,1}(C;L)$'s for all $i\ge r-3$  is roughly the following. Non-zero elements of the $K_{i,1}(C;L)$'s can be seen to correspond to certain subvarieties containing the curve $\varphi_L(C) \subset \P^r$ and defined by quadrics. On the other hand the existence of divisors satisfying $(\Delta_3)$ plays the role of a generality condition which prevents the curve from being contained in such a variety. This simple contradiction works quite efficiently once the curve is projected in $\P^5$, and that's how we prove the theorem. Note that the condition Cliff$(C) \ge 3$ cannot be removed, as easy examples show.

For higher $q$ we have a similar contradiction. But the verification that $(M_q)$ holds once hypotheses similar to those of the theorem are satisfied,  becomes much more involved as $q \ge 4$ and would require a classification of certain classes of varieties that is not available at the moment. 

  It is interesting to note that in Theorem \ref{T:GL1} the existence of secant spaces is related to the exceptions to the validity of $(N_p)$, hence it is a condition not satisfied  in general. On the other hand in Theorem \ref{T:m3} the existence of secant spaces, implied by condition $(\Delta_3)$, is a   condition that is satisfied in general.

 A final note in the case $L=K_C$. The condition Cliff$(C)\ge 2$ already implies the existence of divisors satisfying $(\Delta_2)$. Similarly,
the use of condition $(\Delta_3)$ made by Voisin in \cite{cV88} plays a role in the proof, but is not required for the validity of $(M_3)$: all is required is that
Cliff$(C) \ge 3$; in fact the main difficulty in  \cite{cV88} consists in proving that Cliff$(C) \ge 3$ implies the existence of $D$ satisfying $(\Delta_3)$.  This suggests that, more generally,  Cliff$(C) \ge q$ might imply the existence of divisors $D$ satisfying $(\Delta_q)$ with respect to $K_C$.

The paper is organised as follows. In Section \ref{sec:def} we introduce the main condition $(\Delta_q)$ and study its general properties. In Section \ref{sec:canonical} we specialise to the case of canonical curves. In Section \ref{sec:geom} we relate condition $(\Delta_q)$ to the geometry of the curve in $\mathbb P^r$ and in Section \ref{sec:Koszul} we recall the definition of syzygy schemes and prove Theorem {\ref{T:m3}.

%%%%%%%%%%%%%%%%%%%%%%%%%%%%%%%%%%%%%%%%%%%%%%%%%%%%%%%%%%%%%%%%

\section{The condition $(\Delta_q)$}
\label{sec:def}

\subsection{Secant loci}
For any $n \ge 1$ we denote by   $C_n$  the  $n$-th symmetric product  of $C$ 
 and by $\Xi_n \subset C\times C_n$ the universal divisor. Let 
 \[
 \xymatrix{
 C&C\times C_n\ar[l]_-\pi\ar[d]^{\pi_n} \\
 &C_n}
 \]
 be the projections. For any globally generated line bundle $L$ on $C$, the sheaf on $C_n$: 
 \[
 E_L := \pi_{n*}(\pi^*L\otimes\O_{\Xi_n})
 \]
is locally free of rank $n$ and is called the \emph{secant bundle} of $L$. We have a homomorphism of locally free sheaves   on $C_n$:
\[
\xymatrix{
   \pi_{n*}\pi^*L  \ar@{=}[d]\ar[r]^-{e_{L,n}} & E_L \\
 H^0(L)\otimes\O_{C_n}}
\]
Note that $e_{L,n}$ is generically surjective.

We will denote by $V^k_n(L) \subset C_n$ the closed subscheme defined by the condition
\[
{\rm rank}(e_{L,n}) \le k.
\]
Standard facts about determinantal subschemes, see for example \cite{ACGH}, imply that if non--empty, then $V^k_n(L)$ has dimension $\ge n-(r+1-k)(n-k)$ which is the {\em expected dimension}.

Of special interest are the cases $k=n-1$.  The scheme $V^{n-1}_n(L)$ is supported on the set of $D\in C_n$ which do not impose independent conditions to $L$, and its expected dimension is $2n-r-2$. If $n=r$, we can prove the following:
\begin{lem}\label{L:deltap1}
If $r\ge 4$ then $V_r^{r-1}(L)$ is non--empty and of pure dimension $r-2$.     
%\item[(iii)] If $2p \le r-1$ and $x \in C$ is a general point, then $S_{r-p+1}(L)$ has the expected  dimension if and only if $S_{r-p+1}(L(-x))$ has the expected dimension. 
%\item[(ii)] If $S_{r-q+2}(L)$ has the expected dimension, and $y_1, \dots, y_{r-2q+2}\in C$ are general points then the  projection 
%	$\varphi_{L(-y_1-\cdots-y_{r-2q+2})}(C) \subset \P^{2q-2}$  of  $\varphi_L(C)\subset \P^r$ from $\langle y_1, \dots, y_{r-2q+2} \rangle$ has finitely many 
%	$q$-secant $\P^{q-2}$'s.
\end{lem}

\proof 
%(i) $S_{r-q+2}(L)$ is supported on  the effective divisors $D=x_1+\cdots + x_{r-q+2}$   
%such that the linear span  $\langle D \rangle \subset \P(H^0(L)^\vee)$ has dimension $r-q$. 
%The homomorphism $e_{L,r-q+2}$ has torsion cokernel because $L$ is globally generated.  
%Since the ranks of  $ \pi_{r-q+2*}\pi^*L $ and of $E_L$ are $r+1$ and $r-q+2$ respectively, we have that
%if $S_{r-q+2}(L) \ne \emptyset$ then every component $S'$ of $S_{r-q+2}(L)$ has codim$(S') \le q$.  

Let $\Sigma$ be a non--empty component of $V_r^{r-1}(L)$ with $\mathrm{codim}(\Sigma)\le 1$, i.e. that $\dim(\Sigma) \ge r-1$. Consider the morphisms:
\[
\xymatrix{
C_{r-1} \times C\ar[d]^-{\pi_{r-1}} \ar[r]^-\sigma& C_r \\
C_{r-1}}
\]
Then $\pi_{r-1}(\sigma^{-1}(\Sigma)) = C_{r-1}$. This implies that    if $x_1, \dots,  x_{r-1} \in C$ are general points then  the pencil $|L(-x_1- \cdots -  x_{r-1})|$ has  base points, which is impossible. Therefore $V_r^{r-1}(L)$ has pure dimension $r-2$.  

For the same reason, if $A=x_1+ \cdots + x_{r-2}$ is a general effective divisor of degree $r-2$ then $L(-A)$ is base-point free and not composed with an involution. The plane curve 
$\Gamma:=\varphi_{L(-A)}(C) \subset \P^2$ is singular and birational to $C$. 
Letting $x_{r-1},x_r\in C$ be such that $\varphi_{L(-A)}(x_{r-1}) =\varphi_{L(-A)}(x_r)$ is a singular point of $\Gamma$, the divisor  $x_1+ \cdots +x_{r-2}+ x_{r-1}+x_r$ belongs to $V^{r-1}_r(L)$, which concludes non--emptiness.  
\qed

Let us record the following useful fact which is a direct generalisation of Lemma 1.7, pag. 163 \cite{ACGH}.

\begin{lem}
\label{L:IrredComp}
Assume $q\ge 2$, $r-q+2\ge 4$ and $V^{r-q+1}_{r-q+2}(L)\ne\emptyset$. Then no irreducible component of $V^{r-q+1}_{r-q+2}(L)$ is contained in $V^{r-q}_{r-q+2}(L)$.
\end{lem}

\proof
Let $D=x_1+\cdots +x_{r-q+2}$ be a general element in a component of $V^{r-q+1}_{r-q+2}(L)$. Assume by contradiction that $D\in{V^{r-q}_{r-q+2}(L)}$. Then $\mathrm{dim}\langle D\rangle\le r-q-1$. We may assume that $\langle D\rangle=\langle x_1+\cdots +x_{r-q+1}\rangle$. Then for a general $x\in C$ we have $\mathrm{dim}\langle x_1+\cdots +x_{r-q+1}+x\rangle\le r-q$ and therefore $x_1+\cdots +x_{r-q+1}+x\in V^{r-q+1}_{r-q+2}(L)$. To conclude, note that $x_1+\cdots +x_{r-q+1}+x$ and $D$ belong to the same component of $V^{r-q+1}_{r-q+2}(L)$ and $\mathrm{dim}\langle D\rangle <\mathrm{dim}\langle x_1+\cdots +x_{r-q+1}+x\rangle$ contradicting the generality of~$D$.
\qed

\medskip

A consequence of Lemma \ref{L:IrredComp} is that the locally closed subscheme $S_{r-q+2}(L) \subset C_{r-q+2}$ defined as
\[
S_{r-q+2}(L):= V^{r-q+1}_{r-q+2}(L)\backslash V^{r-q}_{r-q+2}(L)
\]
is dense in any irreducible component of $V^{r-q+1}_{r-q+2}(L)$. In particular, any property which is satisfied by general divisors in any irreducible component of $V^{r-q+1}_{r-q+2}(L)$ is also valid for $S_{r-q+2}(L)$. Note that the expected dimension is $r-2q+2$ in this case. For the particular case $q=2$, Lemma \ref{L:deltap1} shows that the dimension of $S_r(L)$ coincides with the expected dimension $r-2$.

\subsection{Condition $(\Delta_q)$}
We introduce our basic  condition.

\begin{defn}\label{D:deltaq}
Assume that $ r\ge 4$ and let $2\le q \le 1+r/2$. We say that a reduced effective divisor $D = x_1+\cdots + x_{r-q+2}$ on $C$ satisfies condition $(\Delta_q)$ with respect to $L$ if the following conditions are satisfied:
\begin{itemize}
	\item[(a)] $h^0(L(-D)) = q$,
	\item[(b)] $L(-D)$ is base point free,
	\item[(c)] $h^0(L(-D+x_i))=h^0(L(-D))$ for all $i=1,\dots, r-q+2$.
\end{itemize}
\end{defn}

In terms of projective geometry, the conditions defining $(\Delta_q)$ can be rephrased as follows:  
\begin{itemize}
\item[(a)] the linear span $\langle D \rangle \subset \P^r$ is an  $(r-q)$-plane,
\item[(b)] $\langle D \rangle \cap C=\mbox{Supp}(D)$,
\item[(c)] $x_1,\dots, x_{r-q+2}$ are in linearly general position in $\langle D \rangle$ (but not in $\P^r$ of course) i.e. $\langle D-x_i \rangle =\langle D \rangle$ for all $i$.
\end{itemize}

In terms of symmetric products, the conditions defining $(\Delta_q)$ correspond to the following:

\begin{itemize}
\item[(a)] $D\in S_{r-q+2}(L)$,
\item[(b)] $\{D\}+C\subset S_{r-q+3}(L)$,
\item[(c)] $D\not\in \mathrm{Im}\{V^{r-q}_{r-q+1}(L)\times C\to C_{r-q+2}\}$.
\end{itemize}

Note that, from Lemma \ref{L:IrredComp}, a general point in any irreducible component of $V^{r-q+1}_{r-q+2}(L)$ satisfies condition (a). Clearly, divisors $D = x_1+\cdots + x_{r-q+2}$ as in Definition (\ref{D:deltaq}) fill an open subset of  $S_{r-q+2}(L)$.

%In particular, the correspondence $D\mapsto L(-D)$ birationally identifies the components of $V^{r-q+1}_{r-q+2}(L)$ with the components of the locus $W^{q-1}_{d-r+q-2}$, where $d=\mathrm{deg}(L)$.

\begin{prop}
Assume that $L$ is special and embeds $C$ with a $(r-q+2)$-secant $(r-q)$-plane $\langle D\rangle\subset \mathbb P^r$.Then $h^0(\mathcal O_C(D))\le 2$.
\end{prop}

\proof
Assume that $L=K_C(-B)$, and denote by $r_B:=h^0(\mathcal O_C(B))-1\ge 0=h^1(L)-1$. From the Riemann-Roch theorem applied to $L$ we obtain $\mathrm{deg}(B)=r_B-r+g-1$, and hence $\mathrm{deg}(B+D)=g-p+r_B$. From Riemann-Roch applied to $L(-D)$ we obtain $h^0(\mathcal O_C(B+D))=r_B+2$. Since the addition map of divisors $|B|\times |D|\to |B+D|$ is finite on its image, it follows that $\mathrm{dim}|D|\le 1$.

\begin{rmk}\rm\label{R:deltap1}
\begin{itemize}
%\item[(i)] A divisor $D\in C_{r+1}$ satisfying condition $(\Delta_1)$ with respect to $L$ always exists unless $C=\P^1$,  a case  we have excluded from the beginning.   
	 
\item[(i)]  A divisor $D$ satisfies  $(\Delta_q)$ with respect to $K_C$ if and only if $|D|$ is a primitive $g^1_{g-q+1}$.
In particular, $(\Delta_2)$ is equivalent to $|D|$ being  a primitive $g^1_{g-1}$ on $C$, and such a $D$ does not exist if and only if $C$ is  trigonal  or a plane nonsingular quintic (see \cite{GL85}). Note that hyperelliptic curves are excluded automatically by our assumptions if $L=K_C$. We shall treat the canonical case in a separate section.

\item[(ii)]  If $L$ is  non-special of degree   $d=g+r \ge 2g$,  then there is no divisor $D \in C_{r-g+1}$ satisfying condition $(\Delta_{g+1})$ with respect to $L$. In fact this would imply that 
$L(-D)$ is  base-point-free   of degree $(g+r)-(r-g+1)= 2g-1$ and dimension $r-(r-g)=g$, and this is impossible. 
If $g=1$ this means that   no $D\in C_{r}$    
 satisfies $(\Delta_2)$ with respect to $L$: in fact $C\subset \P^r$ has degree $r+1$ and any $r$ distinct points of $C$ are independent.  If $g=0$ we recover Remark (i) above.
\end{itemize}
\end{rmk}

\noindent
\textbf{Terminology:} Assume $L$ to be special and very ample, $h^0(L)=r+1$, and let $2 \le q \le r-1$. It is convenient to introduce the following:
\begin{itemize}
	\item We say that {\em condition $(\Delta_q)$ holds  on a component $V$} of $V^{r-q+1}_{r-q+2}(L)$ if the general element $D\in V$ satisfies $(\Delta_q)$ with respect to $L$. We say that 
	{\em $(\Delta_q)$ holds  on $C$ with respect to $L$}  if it holds on every component of $V^{r-q+1}_{r-q+2}(L)$.
	
	\item We say that $(\Delta_q)$ holds  on $C$ with respect  to $L$ {\em in the strong sense} if it holds and moreover all components of $V^{r-q+1}_{r-q+2}(L)$ have dimension equal to the expected dimension $r-2q+2$.
	A necessary condition for this to happen is that $r\ge 2q-2$.
	\item When we say ``dim$(Z)=d$'' we mean that each irreducible component of $Z$ has dimension $d$.
\end{itemize}

Most of our results are proved only under the assumption that $(\Delta_q)$ holds in the strong sense.

\begin{prop}
\label{prop:DimCond}
Assume that $\mathrm{dim}(V^{r-q+1}_{r-q+2}(L))=r-2q+2$.
Then $(\Delta_q)$ holds on $C$ with respect to $L$ in the strong sense if and only if the following conditions are satisfied:
\begin{enumerate}
\item $\mathrm{dim}(V^{r-q+1}_{r-q+3}(L))\le r-2q+1$,
\item $\mathrm{dim}(V^{r-q}_{r-q+1}(L))=r-2q$.
\end{enumerate}

\end{prop}

\proof
Note that the expected dimension of the locus $V^{r-q+1}_{r-q+3}(L)$ is $r-3q+3\le r-2q+1$.

The proof relies on the observation that {\em any} map defined by addition of divisors is finite on its image.
Assume that $\mathrm{dim
}(V^{r-q+1}_{r-q+3}(L))\le r-2q+1$ and $\mathrm{dim}(V^{r-q}_{r-q+1}(L))=r-2q$.
Let $D\in S_{r-q+2}(L)$ be a general element in an irreducible component. Then, by definition, $h^0( L(-D))=q$, hence condition (a) from Definition \ref{D:deltaq} is satisfied. We prove that $L(-D)$ has no base-points, i.e. condition (b). Suppose that $x$ is a base-point of $L(-D)$, then $D+x\in V^{r-q+1}_{r-q+3}(L)$ and depends on $r-2q+2$ parameters, contradicting the assumption on $\mathrm{dim}(V^{r-q+1}_{r-q+3}(L))$. We have seen that condition (c) is equivalent to $D\not\in \mathrm{Im}\{V^{r-q}_{r-q+1}(L)\times C\to C_{r-q+2}\}$. By the dimensionality assumptions, the image of the addition map cannot fill a dense set of a component of  $S_{r-q+2}(L)$.

Conversely, assume that $(\Delta_q)$ hold on $C$ with respect to $L$ in the strong sense. Suppose that $V^{r-q}_{r-q+1}(L)$ has a component $Z$ with $\mathrm{dim}(Z)\ge r-2q+1$. Then dimensionality hypothesis, the image of the set $Z+C$ inside $V^{r-q+1}_{r-q+2}(L)$ must fill a component, and all its points violate $(\Delta_q)$. If there is a component $Y$ of $V^{r-q+1}_{r-q+3}(L)$ having dimension $\ge r-2q+2$ then for a general element $D'\in Y$ can be written as $D'=D+x$ where, again by dimensionality assumption, $D$ must fill a component of $V^{r-q+1}_{r-q+2}(L)$. From the definition, $D$ fails property (b) of $(\Delta_q)$, contradiction.
%\marginpar{\tiny Prove the converse or remove this implication in the statement.}
\qed
%\marginpar{\tt Excess dimension}

\begin{rmk}
Recall that in the case $q=2$ the dimension of the locus $V_r^{r-1}(L)$ equals the expected dimension $r-2$ (Lemma \ref{L:deltap1})  but it can be reducible:  when $L=K_C$ and $g \ge 6$ this  happens precisely when $C$ is either trigonal or bielliptic (see \cite{mT84}).  In the trigonal case $V_{g-1}^{g-2}(K_C)$ has two components and in both of them $(\Delta_2)$ does not hold. In the bielliptic case $(\Delta_2)$ holds in one component but not in the other.
A characterisation of the pairs $(C,L)$ for which $V_r^{r-1}(L)$ is reducible is unknown to us when $L$ is arbitrary. 
 \end{rmk}

\begin{lem}\label{L:deltap2}
Assume   $r \ge 5$,  $2 \le q \le (r+1)/2$. Assume that $(\Delta_q)$ holds on $C$ with respect to $L$ in the strong sense. Then for every general $x \in C$,  $(\Delta_{q})$ holds on $C$ with respect to $L(-x)$ in the strong sense.
\end{lem}

\proof  As noted before, it suffices to prove the same statement for the locally closed subschemes $S_{r-q+1}$.
Let $x \in C$ be a point such that for each irreducible component of $S_{r-q+2}(L)$ it is not in the support of all divisors of that component and that it is in the support of some divisor in it that satisfies $(\Delta_q)$ with respect to $L$.      We have a diagram of spaces and maps:
\[
\xymatrix{
C_{r-q+1}\times \{x\} \ar[r]^-{\phi} & C_{r-q+2} \\
S_{r-q+1}(L(-x)) \ar[r] \ar[u] & S_{r-q+2}(L)\ar[u]}
\]
 where all the maps are inclusions. Let $\Sigma \subset S_{r-q+1}(L(-x))$ be an irreducible component. Assume that $\dim(\Sigma)\ge r-2q+2$. Then $\phi(\Sigma)$ is a component of $\overline{S_{r-q+2}(L)}$ and all divisors in 
$\phi(\Sigma)$ contain $x$ in their support. This contradicts our assumptions.  The second possibility is that $\dim(\Sigma)= r-2q+1$ and that all divisors $D\in \Sigma$ do not satisfy $(\Delta_{q})$ with respect to $L(-x)$.   Then $\phi(\Sigma) \subset \overline{S_{r-q+2}(L)}$ and all $D+x\in \phi(\Sigma)$ do not satisfy $(\Delta_q)$ with respect to $L$. Since this condition is satisfied for a general choice of $x\in C$ we deduce that there is a component of $\overline{S_{r-q+2}(L)}$ all whose elements do not satisfy $(\Delta_q)$ with respect to $L$, a contradiction. \qed

%\medskip

%\marginpar{\tiny this goes to the canonical section}
%Applying Proposition \ref{prop:canonical} and Lemma \ref{L:deltap2}, we obtain the following existence result:
%
%
%\begin{cor}
%For a general triple $(C,L,D)$, with $L$ special and $D\in V^{r-q+1}_{r-q+2}(L)$, the condition $(\Delta_q)$ is satisfied.
%\end{cor}
%
%
%The meaning of generality in the statement is that $L$ is a general projection of the canonical bundle, and hence the speciality index equals one.
%More precise existence results are proved by Coppens and Martens, or by Farkas, Theorem 0.5~\cite{gF08}.
%

\section{The case $L=K_C$}
\label{sec:canonical}

In this case the notation specialises as follows:
\begin{itemize}
	\item $V_{g-q+1}^{g-q}(K_C)= C^1_{g-q+1}$.
	\item $S_{g-q+1}(K_C)= C^1_{g-q+1}\setminus C^2_{g-q+1}$.
	\item the expected dimension  of $V_{g-q+1}^{g-q}(K_C)$ is $g-2q+1$.
	\item A divisor $D\in C^1_{g-q+1}$ satisfies $(\Delta_q)$ with respect to $K_C$ for some $q \ge 2$ if and only if it defines a \emph{primitive} $g^1_{g-q+1}$ i.e. it is complete, base-point-free and the residual is also base-point-free. 
\end{itemize}

For brevity, when in this section  we say that a {\em condition $(\Delta_q)$ is satisfied} we assume implicitly  ``with respect to $K_C$''.  

The condition $(\Delta_q)$ is well defined in the range $2 \le q \le g-1$. 
When $\left[\frac{g-1}{2}\right] < q \le g-1$  the existence of a $D \in C_{g-q+1}$ satisfying $(\Delta_q)$ is equivalent to the existence of a primitive $g^1_{g-q+1}$ with $g-q+1 < \frac{g+3}{2}$ and therefore $C$ becomes more and more special as $q$ grows because its gonality decreases.
On the other hand, when $2 \le q \le \left[\frac{g-1}{2}\right]$ the condition that there exists $D$ satisfying $(\Delta_q)$ should imply that Cliff$(C) \ge q$ (this is true for $q=2,3$, see the Remark \ref{rmk:Delta&Cliff} below).  
 In this range if this implication is true then the existence of a $D \in C_{g-q+1}$ satisfying $(\Delta_q)$ implies that $C$ is more and more general as $q$ grows.  We are able to clarify this only assuming that $C$ has Clifford dimension one.

\begin{prop}\label{P:symq1}
Assume that $g \ge 2q+1$, $q \ge 2$.   Consider the following conditions:

\begin{itemize}
\item[(i)] The condition   $(\Delta_q)$ holds on $C$ in the strong sense. 
\item[(ii)] $C\subset \P^{g-1}$ is not contained in a $q$-dimensional variety of minimal degree $g-q$.
	 \item[(iii)] For all $1 \le e \le q$ there does not exist a $\bar{D}\in C_{e+1}$ satisfying $(\Delta_{g-e})$.
	\item[(iv)] $\mathrm{gon}(C) \ge q+2$.

\end{itemize}
Then $(i)\Longrightarrow (ii)\Longleftrightarrow (iii) \Longleftrightarrow (iv)$.
\end{prop}

\proof  $(iv) \Longleftrightarrow (iii)$.  $\mathrm{gon}(C) < q+2$ if and only if there exists a primitive $g^1_{e+1}$ for some $1\le e \le q$, and this is equivalent to the existence of  $\bar{D}\in C_{e+1}$ satisfying $(\Delta_{g-e})$.

$(ii) \Longleftrightarrow (iii)$. The existence of a primitive $g^1_{e+1}$ for some $1\le e \le q$ is equivalent to the existence of   $A \in W^1_{q+1}\setminus W^2_{q+1}$, possibly with base points. The union of the linear spans $\langle E \rangle$ as $E \in |A|$ is a $q$-dimensional variety of minimal degree. 

$(i) \Longrightarrow (iii)$. If there exists $\bar{D}\in C^1_{e+1}$ satisfying $(\Delta_{g-e})$, for some $1\le e \le q$ then the locus 
\[
W:=\{\bar{D}+x_1+\cdots+x_{g-(q+e)}:\bar{D}\in C^1_{e+1} \text{\ satisf. $(\Delta_{g-e})$}, x_i\in C\} \subset C_{g-q+1}
\]
 consists of divisors not satisfying $(\Delta_q)$ and has dimension 
\[
\dim(W)\ge g-(q+e)+1 \ge g-2q+1
\]
 Therefore $\overline{W}$ is a component of $C^1_{g-q+1}$, contradicting (i). \qed

\begin{rmk}
The proof of the implication $(i)\Longrightarrow (iii)$ fails  if $g=2q$. In fact   a general curve $C$ of   genus $g=2q$ has a primitive $g^1_{q+1}$  and    $(\Delta_q)$ holds on $C$ in the strong sense.  In this case $V^q_{q+1}(K_C)=C^1_{q+1}$ is reducible in several components of dimension one: their number is given by Castelnuovo's formula (\cite{ACGH}, p. 211).
\end{rmk}

\begin{rmk}
The implication $(ii) \Longrightarrow (i)$ does not hold. In fact if $C$ is a bielliptic curve then $\mathrm{gon}(C)=4$. On the other hand  $C^1_{g-1}$ has two components \cite{mT84}, both having dimension $g-3$, equal to the expected dimension, but $(\Delta_2)$ holds only on one of them.  Therefore in this case the implication holds only in a weak sense.
\end{rmk}

\begin{rmk}
\label{rmk:Delta&Cliff}
If $(\Delta_2)$ holds then Cliff$(C)\ge 2$. This has been proved in \cite{GL85} using Mumford-Martens.  Note that they only assumed that $(\Delta_2)$ holds on \emph{some} component of $C^1_{g-1}$.
The implication  ``$(\Delta_3)$ holds'' $\Longrightarrow$ $\mathrm{gon}(C)\ge 5$ has been considered in \cite{cV88}.  In both cases $q=2,3$ the converse implication 
\[\xymatrix{
\mathrm{Cliff}(C)\ge q \ar@{=>}[r]& (\Delta_q) \ \text{holds on some component of $C^1_{g-q+1}$}}
\] 
has also been proved.  
\end{rmk}
%This is what we still do not have in general, both for $q \ge 4$ and $L=K_C$ and for $q \ge 2$ and $L$ any special line bundle.

\begin{rmk}
Assume $g$ odd. On a  general curve $C$ of Clifford dimension one there is a $D\in C_{\frac{g+3}{2}}$ satisfying $(\Delta_{\frac{g-1}{2}})$.
 The reason is that $C$ has gonality  $\frac{g+3}{2}$, and a pencil computing its gonality is necessarily primitive. Therefore a divisor $D$ in the pencil satisfies $(\Delta_{\frac{g-1}{2}})$.
\end{rmk}

%\begin{thm}\label{T:deltabyind}
%Assume that $g \ge 2q+2$, $q \ge 3$, and that  $\mathrm{dim}(C^2_{g-q+2})\le....$. If   $(\Delta_q)$ holds on $C$ in the strong sense then  $(\Delta_{q-1})$ also holds on $C$ in the strong sense. 
%\end{thm}

%\proof   By hypothesis  $C^1_{g-q+1}$ has pure dimension $g-2q+1$ or, equivalently, 
%\[
%\dim(W^1_{g-q+1})=g-2q = \rho(g,1,g-q+1)
%\]
%  Then 
%\[
%\dim(W^1_{g-q+2})\ge \rho(g,1,g-q+2)=g-2q+2 = \dim(W^1_{g-q+1})+2
%\]
%By the excess linear series bound we have equality on each component of $W^1_{g-q+2}$, and therefore all components of $C^1_{g-q+2}$ have dimension equal to the expected dimension $g-2q+3$.   
%Moreover 
%\[
%\dim(W^1_{g-q+1}+C) = g-2q+1
%\]
%and therefore the general element of each component of $C^1_{g-q+2}$ satisfies conditions (a) and (c) of Definition \ref{D:deltaq}.
%
%For condition (b) we use the dictionary after Definition \ref{D:deltaq}.
%
%\qed
 
Using "excess linear series", in the case $L=K_C$, Proposition \ref{prop:DimCond} implies:

\begin{prop}
\label{prop:canonical}
Let $C$ be a curve of genus $g\ge 2q+2$ such that the dimension of the locus $W^1_{g-q}(C)$ equals the expected dimension $g-2q-2$ and $\mathrm{dim}(W^2_{g-q+2}(C))\le g-2q-2$. Then  $(\Delta_q)$ holds on $C$ in the strong sense. 
\end{prop}

\begin{rmk}\rm
If the curve $C$ is of gonality $(q+1)$ or less, then the hypotheses before are not satisfied. Indeed, if $A$ is a $g^1_{q+1}$, then $W^1_{g-q}(C)$ contains the variety $\{A\}+W_{g-2q-1}(C)$ which is of dimension $g-2q-1$. 

If the curve $C$ is instead of gonality $(q+2)$, then the hypothesis that $$\mathrm{dim}(W^1_{g-q}(C))=\rho(g,1,g-q)=g-2q-2$$ coincides with the {\em linear growth condition} on the dimension of Brill-Noether loci, from \cite{AproduMRL}. It was proved in \cite{AproduMRL} that this condition implies Green's conjecture, i.e. condition $(M_q)$.
\end{rmk}

\begin{rmk}\rm
If $q=2$, and $C$ is neither trigonal, bielliptic nor plane quintic, the hypotheses of Proposition \ref{prop:canonical} are satisfied. Indeed, if one of the two fails, then we obtain a contradiction with Mumford--Martens' dimension theorem. Likewise, for $q=3$ the failure of the hypotheses contradicts Keem's dimension theorem, \cite{cV88}, Proposition II.0.
\end{rmk}

\begin{cor}\label{T:deltabyind}
Assume that $g \ge 2q+2$, $q \ge 2$, and that  $\mathrm{dim}(C^2_{g-q+2})\le g-2q$. If   $(\Delta_{q+1})$ holds on $C$ in the strong sense then  $(\Delta_{q})$ also holds on $C$ in the strong sense. 
\end{cor}

Applying Proposition \ref{prop:canonical} and Lemma \ref{L:deltap2}, we obtain the following existence result:

\begin{cor}
For a general triple $(C,L,D)$, with $L$ special and $D\in V^{r-q+1}_{r-q+2}(L)$, the condition $(\Delta_q)$ is satisfied.
\end{cor}

The meaning of generality in the statement is that $L$ is a general projection of the canonical bundle, and hence the speciality index equals one.
More precise existence results are proved by Coppens and Martens, and by Farkas, see Theorem 0.5~\cite{gF08} and the references therein.

\section{Condition $(\Delta_q)$ and geometry}
\label{sec:geom}

\begin{prop}\label{P:mp1} 
Assume that  $r \ge \mathrm{max}\{4,2q-1\}$.  Suppose that  $L$ is special and condition $(\Delta_q)$ holds on $C$ with respect to $L$ in the strong sense. Then  $\varphi_L(C)\subset \P^r$ is not contained in a $q$-dimensional variety of minimal degree $(r-q+1)$ unless: $r=2q-1$ and $C$ has a base-point-free $g^1_{q+1}$.% and $V^{r-q+1}_{q+1}(L)$ has a rational component.
 \end{prop}

\proof      Assume that $r\ge 2q$. We note that $C$ has no $g^1_{q+1}$. Indeed, if we have a $g^1_{q+1}$, then $A+C_{r-2q+1}$, with $A\in|g^1_{q+1}|$, fill up a component of $V^{r-q+1}_{r-q+2}(L)$ and any element of this locus fails condition (c) from the definition of $(\Delta_q)$.

Assume by contradiction that $\varphi_L(C)\subset X$,  a $q$-dimensional variety of minimal degree $r-q+1$. Then $X$ is ruled by a one-dimensional family of $(q-1)$-planes. Let $\Lambda$ be a general such $(q-1)$-plane, and let $E=\Lambda \cap   \varphi_L(C)$ and $n=\deg(\Lambda \cap   \varphi_L(C))$. Then $n \ge q+2$ by what we have just shown. Decompose $E= A+B$ with $\deg(A)=q+1$. 
%We may assume that $\langle A \rangle =\Lambda$. 
Let $D=A+y_1+\cdots+y_{r-2q+1}$ with the $y_i$'s general points of $C$. Then $D \in V^{r-q+1}_{r-q+2}(L)$, but it does not satisfy $(\Delta_q)$. On the other hand the divisor $D$ depends on $1+(r-2q+1)=r-2q+2$ parameters. Therefore it is a general point of a component of $V^{r-q+1}_{r-q+2}(L)$, a contradiction.  
		
		In the case $r=2q-1$, the only possibility for $C$ to be on a variety of minimal degree is that $C$ has a base-point-free $g^1_{q+1}$ and, in this case, $S_{q+1}(L)$ will have a rational component.
		The case when  $X$ is a cone over the Veronese surface can be treated similarly, by general projection to $\mathbb P^{2q-1}$ using Lemma \ref{L:deltap2}. 		
		\qed

%In particular, Proposition \ref{P:mp1} applies when a general $D$ in any component of $S_{r-q+2}(L)$ satisfies $(\Delta_q)$ with respect to $L$.

Note that if $C$ is contained in an $e$--dimensional variety of minimal degree $(r-e+1)$ with $e\le q$ then it is contained also in a $q$--dimensional variety of minimal degree $(r-q+1)$  \cite{jH81}.

\medskip

 As we will see, the validity of property $(M_3)$ is tightly connected with properties of  surfaces of low degree in $\P^5$. As an   illustration of the geometric content of Definition \ref{D:deltaq} we study surfaces of degree $n \le 6$.

\begin{prop}\label{P:cast}
 Assume that $r=5$  and that $L$ is special and $(\Delta_3)$ holds on $C$ with respect to $L$ in the strong sense.  
Then $\varphi_L(C)\subset \P^5$  is not contained in a nonsingular surface of degree $\le 6$ unless it has a $g^1_4$.
\end{prop}

\proof  Assume that $C \subset S$, nonsingular of degree $n \le 6$.
Consider the case $n=6$.  The possibilities for a nonsingular surface of degree 6 in $\P^5$ are described in \cite{pI84} and are the following: (i) An elliptic scroll with sectional genus $g = 1$ and $e=0$; (ii) A Castelnuovo surface, with sectional genus $g = 2$, defined by the embedding  in $\P^5$ of the blow-up $X=Bl_{p_1,\dots,p_6,q}(\P^2)$ of $\P^2$ at seven general points via the very ample linear system 
$|\mathcal{L}| = |4H - E_1 - E_2 -\cdots - E_6-2A|$  (with obvious notation) corresponding to the system of plane quartics passing simply through $p_1, \dots, p_6$ and doubly through $q$. 

In case (i) let $\ell\subset S$ be a general line of the ruling and let $k=\mathrm{deg}(\mathcal O_C(\ell))$. Then $k \ge 2$ and if $k = 2$ then $C$ is bielliptic, so it has a $g^1_4$. If $k \ge 3$ then adding a general $p\in C$ to a subdivisor of degree 3 of $\mathcal{O}_C(\ell)$ we obtain an element of $S_4(L)$ which does not satisfy $(\Delta_3)$ and it depends on two parameters, a contradiction.

%\marginpar{\tiny decide whether is A or 2A}
\noindent
In case (ii) the system $|H-A|$ is a pencil of conics on the surface $S$. The divisors $D\in |\mathcal{O}_C(H-A)|$  have degree say $m\ge 3$ and dim$|D|\ge 1$. 
If $m \le 4$ then $C$ has a $g^1_4$. Otherwise the divisors $D$ contain subdivisors of degree 4 contradicting the other conditions. 

If $n=5$ then $S$ is a Del Pezzo surface. Let $|\gamma|$ be a pencil of conics on $S$ and let $N=\mathcal{O}_C(\gamma)$. Then $N$ gives a $g^1_4$ or contradicts $(\Delta_3)$ according to the possibilities $\mathrm{deg}(N) \le 4$ or $\mathrm{deg}(N)\ge 5$.

If $n=4$ the conclusion follows from Proposition \ref{P:mp1}. 

\qed

\section{Condition $(\Delta_3)$ and Koszul cohomology}
\label{sec:Koszul}

In this section we briefly recall the relation between Koszul cohomology and vector bundles, as well as the definition of syzygy schemes.

Consider $X$ a smooth projective variety and let $L$ be a globally generated line bundle on $X$. We let
\[
\varphi_L: X \to \P(H^0(L)^\vee)  \cong \P^r,  \quad  r+1 = h^0(L)
\]
be the morphism defined
by $L$.

			We have an exact sequence
\begin{equation}\label{E:prel1}
\xymatrix{ 0 \ar[r]& M_L \ar[r]& H^0(L) \otimes \O_X \ar[r]& L
\ar[r]& 0}
\end{equation}
where $M_L=\varphi^*(\Omega_{\mathbb P^r}(1))$ is locally free of rank $r$.  If $r=1$, i.e. if $|L|$ is a base point free pencil, then $M_L = L^{-1}$. 
Taking $n$-th exterior power ($1 \le n \le r$) we obtain the exact sequence:
\begin{equation}\label{E:wed1}
	\xymatrix{
	0 \ar[r]& \bigwedge^n M_L \ar[r] & \bigwedge^n H^0(L) \otimes\O_X \ar[r] & \bigwedge^{n-1}M_L \otimes L \ar[r]&0}
\end{equation}

For any coherent sheaf $\mathcal F$ on $X$, twisting the sequence above with $\mathcal F$, with powers of $L$ and taking global sections, we obtain isomorphisms
\[
K_{n,m}(X,\mathcal F;L)\cong \mathrm{Coker}\left\{\wedge^{n+1}H^0(L)\otimes H^0(\mathcal F\otimes L^{m-1})\to H^0(\wedge^nM_L\otimes \mathcal F\otimes L^m)\right\}.
\]

The syzygy schemes were introduced and studied in \cite{mG84}, \cite{sE94}. The idea behind the definition of syzygy schemes is that one reason for which a linearly normal curve $C$ in $\mathbb P^r$ has some nonvanishing $K_{n,1}$  is that $C$ lies on a variety of special type. The varieties under question are cut out by quadrics, more precisely by the quadrics involved in syzygies.

The general set-up is the following. Let $C$ be a smooth curve, $L$ a globally generated (preferably very ample) line bundle on $C$ and denote $V=H^0(L)$.
Start with the short exact sequence of sheaves on the projective space
\[
0\to \mathcal I_C\to \mathcal O_{\mathbb P^r}\to \mathcal O_C\to 0.
\]
Note that for any $n$ and $m$ we have $K_{n,m}(\mathbb P^r,\mathcal O_C;\mathcal O_{\mathbb P^r}(1))\cong K_{n,m}(C;L)$.

Taking Koszul cohomology with respect to $\mathcal O_{\mathbb P^r}(1)$, and using the vanishing of Koszul cohomology on the projective space, we obtain isomorphisms
\[
K_{n,m}(C;L)\cong K_{n-1,m+1}(\mathbb P^r,\mathcal I_C;\mathcal O_{\mathbb P^r}(1)),
\]
for any $n$ and $m$ except for the cases $(n,m)=(0,0)$ or $(n,m)=(1,-1)$. 
On the other hand, from the general description of mixed Koszul cohomology, we know that
\[
K_{n-1,m+1}(\mathbb P^r,\mathcal I_C;\mathcal O_{\mathbb P^r}(1))\cong\mathrm{Coker}\left\{\wedge^nV\otimes H^0(\mathcal I_C(m))\to H^0(\Omega^{n-1}_{\mathbb P^r}(n+m)\otimes\mathcal I_C)\right\}.
\]

Observe that for the case $m=1$, we have $H^0(\mathcal I_C(m))=0$, and hence we obtain an isomorphism
\[
K_{n,1}(C;L)\cong H^0(\Omega^{n-1}_{\mathbb P^r}(n+1)\otimes\mathcal I_C),
\]
in particular, any nonzero Koszul cohomology class $\alpha\in K_{n,1}(C;L)$ corresponds to a section in $H^0(\Omega^{n-1}_{\mathbb P^r}(n+1))$ vanishing along $C$. The zero-scheme of this section is called the {\em syzygy scheme} associated to $\alpha$, and is denoted by $\mathrm{Syz}(\alpha)$. Note that a syzygy scheme is cut out by quadrics, as the sheaf $\Omega^{n-1}_{\mathbb P^r}(n+1)$ is a subsheaf of $\wedge^{n-1}V\otimes \mathcal O_{\mathbb P^r}(2)$.
The scheme-theoretic intersection of all the syzygy schemes is denoted by $\mathrm{Syz}_n(C)$. It contains $C$ and is cut out by quadrics as well.

We record next the following remarkable two classification results concerning syzygy schemes due to Green and Ehbauer.

\begin{thm}[Green's $K_{p,1}$]
If $K_{r-1,1}(C,L)\ne 0$ then $C$ is a rational normal curve and $\mathrm{Syz}_{r-1}(C)=C$. If $C$ is of degree $\ge r+2$ and $K_{r-2,1}(C,L)\ne 0$ then $\mathrm{Syz}_{r-2}(C)$ is a surface of minimal degree $(r-1)$.
\end{thm}

\begin{thm}[Ehbauer]
\label{thm:Ehbauer}
If $C$ has degree $\ge r+13$ and $K_{r-3,1}(C,L)\ne 0$ then $\mathrm{Syz}_{r-3}(C)$ is either a surface of minimal degree $(r-1)$ or a surface of degree $r$ or a 3-fold of minimal degree $(r-2)$.
\end{thm}

We recall the following:

\begin{defn}
\label{def:Mq}
The line bundle $L$ \emph{has property} $(M_q)$ if $K_{n,1}(C;L)=0$ for all $n\ge r-q$.
\end{defn}

We prove:

\begin{thm}{\label{T:m3}}
Assume $g\ge 14$, $r \ge 5$, that $L$ is very ample and special of degree $\ge r+13$,  and  that   $(\Delta_3)$ holds on $C$
with respect to $L$ in the strong sense. Then $L$ satisfies $(M_3)$ unless $\mathrm{gon}(C) \le 4$.
\end{thm}

\proof
Applying Ehbauer's characterisation of syzygy schemes, if $L$ fails property $(M_3)$ then $C$ lies either on a surface of minimal degree or on a 3-fold of minimal degree or on a surface of degree $r$. The first two cases are excluded by Proposition \ref{P:mp1}. Projecting generically to $\mathbb P^5$ and applying lemma \ref{L:deltap2} and Proposition \ref{P:cast}, we see that $C$ cannot lie on a smooth surface of degree $5$. If it lies on a singular surface of degree $5$ in $\mathbb P^5$ then projecting from a singular point, the curve in $\mathbb P^4$ lies on a surface of minimal degree. In particular, since the curve is of gonality $\ge 5$, the image of $C$ in $\mathbb P^4$ has a $k$-secant line for a $k\ge 5$, and hence the image of $C$ in $\mathbb P^5$ has a one-dimensional family of $k$-secant $2$-planes with $k\ge 5$ which contradicts the assumptions. \qed

\begin{rmk}\rm
The same argument together with Green's $K_{p,1}$-Theorem gives a similar statement for the weaker property $(M_2)$. 
%A slightly stronger theorem, proved in \cite{AS78}, states that $(M_2)$ holds for a special $L$ if $\varphi_L(C)$ is not contained in a surface of minimal degree.
\end{rmk}

\end{document}